# Gauss quadrature for integrals and sums

A. D. Alhaidari

*Saudi Center for Theoretical Physics, P.O. Box 32741, Jeddah 21438, Saudi Arabia*

**Abstract:** Gauss quadrature integral approximation is extended to include integrals with a measure consisting of continuous as well as discrete components. That is, we give an approximation for the integral of a function plus its sum over a discrete weighted set.

**Keywords**: Gauss quadrature, integral approximation, continuous measure, discrete measure, mixed measure, orthogonal polynomials, recursion relation.

## 1. Introduction

Approximating integrals, especially intractable ones, by finite sums is as old of a problem as calculus itself. Most approximation procedures follow the structure of the original classic discretized approximation of the integral of a real function $f(x)$ defined in $x_- \leq x \leq x_+$ as

$$\int_{x_-}^{x_+} f(x)\,dx \cong \sum_{n=0}^{N-1} \Delta_n f(x_n), \qquad (1)$$

for some large enough integer $N$. The accuracy of the approximation and its convergence have always been at the core of the development of such approximation schemes. One of the most prominent techniques used in these approximations is the celebrated Gauss quadrature integral approximation [1-3]. It exploits the orthogonality and recursion relation of orthogonal polynomials. The quadrature has also been extended to the approximation of infinite sums (see, for example, [4] and [5] and references therein). In this work, we combine the two. That is, we propose an extended Gauss quadrature associated with integral measures that consist of discrete as well as continuous components. In Section 2, we summarize the conventional Gauss quadrature (GQ) integral approximation where the integration measure is purely continuous. In Section 3, we summarize GQ for approximating infinite sums by finite ones. Finally, in Section 4, we present our main findings where we combine the two.

## 2. The conventional GQ: Purely continuous measure

Let $f(x)$ be an integrable real valued function with respect to the positive measure $d\mu(x) = \rho(x)dx$ in a real vector space, which is spanned by the complete set of orthonormal basis $\{p_n(x)\}_{n=0}^{\infty}$, where $x \in [x_-, x_+] \subset \Re$. Orthogonality is defined as

$$\int_{x_-}^{x_+} \rho(x) p_n(x) p_m(x)\,dx = \delta_{n,m}. \qquad (2)$$

This suggests that we could take $p_n(x)$ to be an orthogonal polynomial of degree $n$ in $x$ with pure continuous spectrum and where $\rho(x)$ is the associated positive weight function [6-8]. The spectral theorem (a.k.a. Favard's theorem) dictates that such polynomials satisfy the following symmetric three-term recursion relation

$$x\,p_n(x) = a_n p_n(x) + b_{n-1} p_{n-1}(x) + b_n p_{n+1}(x), \tag{3}$$

for $n = 1, 2, 3, \ldots$ and where the "recursion coefficients" $\{a_n, b_n\}$ are real constants such that $b_n \neq 0$ for all $n \geq 0$. This recursion gives all the polynomials of any degree starting with the two seed values $p_0(x) = 1$ and $p_1(x) = (x - a_0)/b_0$ (i.e., $b_{-1} \equiv 0$). Associated with this space is an infinite dimensional real tridiagonal symmetric matrix (a Jacobi matrix) whose elements are

$$J_{n,m} = a_n \delta_{n,m} + b_{n-1} \delta_{n,m+1} + b_n \delta_{n,m-1} = \int_{x_-}^{x_+} \rho(x) p_n(x) \, x \, p_m(x) \, dx, \tag{4a}$$

$$J = \begin{pmatrix} a_0 & b_0 & & & \\ b_0 & a_1 & b_1 & & \\ & b_1 & a_2 & b_2 & \\ & & \times & \times & \times \\ & & & \times & \times \end{pmatrix}. \tag{4b}$$

One can easily show that by iterating (4a) for $x$, $x^2$, $x^3$,... and using the three-term recursion relation (3) repeatedly, we obtain

$$\int_{x_-}^{x_+} \rho(x) p_n(x) \, x^k \, p_m(x) \, dx = \left(J^k\right)_{n,m}. \tag{4c}$$

For numerical computations, however, the space is truncated to a finite $N$-dimensional subspace spanned by $\{p_n(x)\}_{n=0}^{N-1}$. The tridiagonal matrix (4) becomes a finite $N \times N$ matrix $J$. The real $N$ distinct eigenvalues of $J$, which we designate as the set $\{\varepsilon_n\}_{n=0}^{N-1}$, are the zeros of the polynomial $p_N(x)$ [i.e. $p_N(\varepsilon_n) = 0$]. Let $\{\Lambda_{m,n}\}_{m=0}^{N-1}$ be the normalized eigenvector of $J$ associated with the eigenvalue $\varepsilon_n$. In this setting, Gauss quadrature integral approximation states that [1-3]

$$\int_{x_-}^{x_+} \rho(x) f(x) \, dx \cong \sum_{n=0}^{N-1} \omega_n f(\varepsilon_n), \tag{5}$$

where the "numerical weights" could be evaluated as $\omega_n = \Lambda_{0,n}^2$. Due to the lower numerical cost in computing matrix eigenvalues instead of eigenvectors, we can also write these numerical weights in terms of $\{\varepsilon_n\}_{n=0}^{N-1}$ and another set of eigenvalues $\{\hat{\varepsilon}_m\}_{m=0}^{N-2}$ as

$$\omega_n = \frac{\displaystyle\prod_{m=0}^{N-2} (\varepsilon_n - \hat{\varepsilon}_m)}{\displaystyle\prod_{\substack{k=0 \\ k \neq n}}^{N-1} (\varepsilon_n - \varepsilon_k)}, \tag{6}$$

where $\{\hat{\varepsilon}_m\}_{m=0}^{N-2}$ is the set of eigenvalues of a submatrix of $J$ obtained by deleting the first (zeroth) row and first column. If sorted, these eigenvalues interlace as $\varepsilon_0 < \hat{\varepsilon}_0 < \varepsilon_1 < \hat{\varepsilon}_1 < \varepsilon_2 \ldots \hat{\varepsilon}_{N-2} < \varepsilon_{N-1}$. The integral approximation (5) becomes exact if $f(x)$ is a polynomial in $x$ of a degree less than or equal to $2N-1$. If, instead of the integral (5), we have $\int_{x_-}^{x_+} f(x)dx$ then this could be approximated as follows

$$\int_{x_-}^{x_+} f(x)dx = \int_{x_-}^{x_+} \rho(x)\frac{f(x)}{\rho(x)}dx \cong \sum_{n=0}^{N-1} \omega_n \frac{f(\varepsilon_n)}{\rho(\varepsilon_n)} := \sum_{n=0}^{N-1} \tilde{\omega}_n f(\varepsilon_n), \tag{7}$$

where $\{\tilde{\omega}_n\}$ are named the "derivative weights" and $\tilde{\omega}_n = \omega_n/\rho(\varepsilon_n)$.

Very often, one encounters integrals that represent matrix elements of functions in the basis set $\{\phi_n(x)\} = \{\sqrt{\rho(x)}\, p_n(x)\}$ (for example, the matrix elements of a potential function in quantum mechanics) of the form

$$f_{n,m} = \int_{x_-}^{x_+} \phi_n(x)f(x)\phi_m(x)dx = \int_{x_-}^{x_+} \rho(x)p_n(x)f(x)p_m(x)dx. \tag{8}$$

Using Eq. (5), this integral is approximated as follows

$$f_{n,m} \cong \sum_{k=0}^{N-1} \omega_k p_n(\varepsilon_k)f(\varepsilon_k)p_m(\varepsilon_k). \tag{9}$$

Now, it is well established that $p_n(\varepsilon_k) = \Lambda_{n,k}/\Lambda_{0,k}$ for $n,k = 0,1,2,..,N-1$. Therefore, with $\omega_k = \Lambda_{0,k}^2$ this could be rewritten in matrix form as

$$f_{n,m} \cong \sum_{k=0}^{N-1} \Lambda_{n,k} f(\varepsilon_k) \Lambda_{m,k} = \left(\Lambda F \Lambda^{\mathrm{T}}\right)_{n,m}, \tag{10}$$

where $F$ is a diagonal matrix with elements $F_{k,j} = f(\varepsilon_k)\delta_{k,j}$. Therefore, to obtain an approximate evaluation of integrals using Gauss quadrature associated with orthogonal polynomials satisfying (2) and (3), one needs only the tridiagonal symmetric matrix $J$, which is constructed using the recursion coefficients $\{a_n, b_n\}$, and possibly the weight function $\rho(x)$ for integrals of the type (7).

## 3. GQ for a purely discrete measure

If instead of the integral $\int_{x_-}^{x_+} \rho(x)f(x)dx$, we were tasked with approximating the infinite weighted sum $\sum_{k=0}^{\infty} \xi_k f(x_k)$, for a given discrete set $\{\xi_k, x_k\}$, by a finite sum similar to that of Eq. (5). Then, we could approach the solution to the problem as follows. We can convert the integral $\int_{x_-}^{x_+} \rho(x)f(x)dx$ into the sum $\sum_{k=0}^{\infty} \xi_k f(x_k)$ by discretizing the weight function and writing $\rho(x) = \chi(x)\sum_{k=0}^{\infty} \delta(x-x_k)$ with $\xi_k = \chi(x_k)$. Performing this same weight function discretization in the orthogonality (2) turns it into

$$\sum_{k=0}^{\infty} \xi_k p_n(x_k)p_m(x_k) = \delta_{n,m}. \tag{11}$$

There are many discrete polynomials that satisfy this kind of orthogonality. The Appendix lists relevant properties of few of these discrete hypergeometric orthogonal polynomials that belong to the Askey scheme [9]. Some of these polynomials constitute an infinite discrete sequence whereas others form finite sequences. Examples of the former are the Charlier and Meixner polynomials, whereas the Krawtchouk polynomial is an example of the latter [9]. All these polynomials do satisfy symmetric three-term recursion relations similar to (3) that read

$$x_k p_n(x_k) = a_n p_n(x_k) + b_{n-1} p_{n-1}(x_k) + b_n p_{n+1}(x_k).  \tag{12}$$

Consequently, we can write the following approximation

$$\sum_{k=0}^{M} \xi_k f(x_k) \cong \sum_{n=0}^{N-1} \omega_n f(\varepsilon_n),  \tag{13}$$

where $M$ can be a finite or an infinite integer depending on the spectrum size of the polynomials $\{p_n(x_k)\}$. However, typically $M \gg N$ otherwise the approximation is not called for. Moreover, the set $\{\varepsilon_n, \omega_k\}_{n=0}^{N-1}$ in (13) are calculated using the tridiagonal symmetric matrix $J$ obtained from (12) exactly as in Section 2 and shown by Eq. (4b). Additionally, the approximation (13) becomes exact if $f(x)$ is a polynomial in $x$ of a degree less than or equal to $2N-1$. Here too, if instead of the sum (13) we have $\sum_{k=0}^{M} f(x_k)$, then we can write

$$\sum_{k=0}^{M} f(x_k) \cong \sum_{n=0}^{N-1} \tilde{\omega}_n f(\varepsilon_n),  \tag{14}$$

where the derivative weights are calculated as $\tilde{\omega}_n = \omega_n / \chi(\varepsilon_n)$. As an example, we consider the function $f(x) = r^x / \Gamma(x+1)$ with $(x, r) > 0$ and obtain the approximation of its infinite sum (14) using GQ associated with the Charlier and Meixner polynomials. Since $x_k = k = 0, 1, 2, ...$ for both polynomials, then the infinite sum has the exact value $e^r$. Table 1 demonstrates the validity of formula (14) and illustrates the relative accuracy and convergence associated with the two polynomials for this particular function. In another example, we consider the finite sum (14) for the function $f(x) = (x+1)r^{x+1} / \Gamma(x+r+2)$ with $M = 100$ and employ GQ associated with the finite Krawtchouk polynomial. Since $x_k = k$, then this finite sum has the following exact closed form (see formula 3, Section 1.3.9 in Ref. [10])

$$\sum_{k=0}^{M} \frac{(k+1)r^{k+1}}{\Gamma(k+r+2)} = \frac{1}{\Gamma(r)} - \frac{r^{M+2}}{\Gamma(M+r+2)}.  \tag{15}$$

Table 2 demonstrates the validity of formula (14) and illustrates its accuracy as well as its convergence as we vary $N$ and the Krawtchouk polynomial parameter $\gamma$.

## 4. GQ for a measure with a mix of continuous and discrete spectra

If the integral measure $d\mu(x) = \rho(x)dx$ contains a mix of discrete as well as continuous spectra, then the weight function could be written as follows

$$\rho(x) = \sigma(x) + \chi(x) \sum_{k=0}^{M} \delta(x - x_k),  \tag{16}$$

where $\sigma(x)$ is the continuous component of the weight function and $\chi(x_k) := \xi_k$ is the discrete component. The discrete spectrum (a.k.a. mass points or bound states) is the set $\{x_k\}_{k=0}^M$ of size $M+1$, which could be finite or infinite. The orthogonality (2) becomes [6-8]

$$\int_{x_-}^{x_+} \sigma(x) p_n(x) p_m(x) dx + \sum_{k=0}^{M} \xi_k p_n(x_k) p_m(x_k) = \delta_{n,m}. \tag{17}$$

These polynomials do also satisfy a symmetric three-term recursion relation similar to (3) where the associated infinite tridiagonal symmetric matrix $J$ is endowed with both continuous as well as discrete eigenvalues. Appendix A lists the relevant properties of two of these hypergeometric orthogonal polynomials that belong to the Askey scheme [9].

If the function $f(x)$ is defined over $x \in [x_-, x_+] \cup \{x_k\}_{k=0}^M$ then Gauss quadrature integral approximation (5) could be extended for this scenario to read

$$\int_{x_-}^{x_+} \sigma(x) f(x) dx + \sum_{k=0}^{M} \xi_k f(x_k) \cong \sum_{n=0}^{N-1} \omega_n f(\varepsilon_n), \tag{18}$$

where the set $\{\varepsilon_n, \omega_k\}_{n=0}^{N-1}$ in (18) are calculated using the tridiagonal symmetric matrix $J$ exactly as in Section 2. Here too, the integral approximation (18) becomes exact if $f(x)$ is a polynomial in $x$ of a degree less than or equal to $2N-1$. If the sum in (18) is finite such that $M$ is of the same order as $N$ (i.e., $M \lesssim N$), then one could think of (18) as an integral approximation by the difference between two finite sums as follows

$$\int_{x_-}^{x_+} \sigma(x) f(x) dx \cong \sum_{n=0}^{N-1} \omega_n f(\varepsilon_n) - \sum_{k=0}^{M} \xi_k f(x_k). \tag{19}$$

As an example, we consider the function $f(x) = x^m e^{-x/2}$ where $m$ is a positive integer and obtain an approximation of its integral-sum (18) using GQ associated with the continuous dual Hahn polynomial $S_n^\mu(x^2; \alpha, \alpha)$ with a mixed spectrum (i.e., $\mu < 0$). Since the polynomial argument is not $x$ but $x^2$ and $x_k^2 = -(k+\mu)^2$, then the GQ integral formula (18) must be revised to read

$$\int_0^\infty \sigma(x) f(x^2) dx + \sum_{k=0}^{M} \xi_k f(x_k^2) \cong \sum_{n=0}^{N-1} \omega_n f(\varepsilon_n). \tag{20}$$

Table 3 demonstrates the validity this formula and illustrates its convergence for a given $\mu$ and several values of $\alpha$. The GQ result is compared to the exact result $[f(J)]_{00}$ [11].

We plan to follow this work by another in which we attempt at approximating the following integral-sum combination that replaces (18)

$$\int_{x_-}^{x_+} f(x) dx + \sum_{k=0}^{M} f(x_k) \cong \sum_{n=0}^{N-1} \tilde{\omega}_n f(\varepsilon_n), \tag{21}$$

where we determine the derivative weights $\tilde{\omega}_n$ in terms of $\omega_n$, $\sigma(\varepsilon_n)$ and $\chi(\varepsilon_n)$.

# Appendix: Relevant orthogonal polynomials

In this Appendix, we give the most essential properties of orthogonal polynomials that are relevant to our work. All these polynomials belong to the hypergeometric class in the Askey scheme [9]. We consider the orthonormal version of these polynomials. That is, where the corresponding three-term recursion relation is symmetric as given by (3) and the orthogonality is a pure $\delta_{n,m}$ similar to those given by Eqs. (2), (11) and (17). Therefore, such polynomials are well-defined once we give the recursion coefficients $\{a_n, b_n\}$ and the weight function $\rho(x)$, $\{\xi_k\}$, or the pair $\{\sigma(x), \xi_k\}$ corresponding to continuous, discrete, or mixed spectra, respectively. Moreover, we always take the initial polynomials as $p_0(x) = 1$ and $p_1(x) = (x - a_0)/b_0$.

## A.1 Discrete weight function

In this section of the Appendix, we give the essential properties of three of the discrete polynomials. The first two form infinite sequences whereas the third is a finite sequence. However, all have $x_k = k$.

*The Charlier polynomial:*

$$C_n^\mu(k) = \sqrt{\frac{\mu^n}{n!}} \,_2F_0\left(\begin{matrix}-n,-k\\ \_\end{matrix}\bigg| -1/\mu\right), \tag{A1}$$

where $\mu > 0$ and $k, n = 0, 1, 2, \ldots$. The recursion coefficients are

$$a_n = n + \mu, \qquad b_n = -\sqrt{\mu(n+1)}. \tag{A2}$$

The discrete weight function reads as follows

$$\xi_k = \mu^k / e^\mu k!. \tag{A3}$$

*The Meixner polynomial:*

$$M_n^\mu(k;\beta) = \sqrt{\frac{(2\mu)_n}{n!}} \beta^{n/2} \,_2F_1\left(\begin{matrix}-n,-k\\ 2\mu\end{matrix}\bigg| 1-\beta^{-1}\right), \tag{A4}$$

where $\mu > 0$, $1 > \beta > 0$ and $k, n = 0, 1, 2, \ldots$. Moreover, $(a)_n = a(a+1)(a+2)\ldots(a+n-1) = \frac{\Gamma(n+a)}{\Gamma(a)}$ is the Pochhammer symbol (the shifted factorial). The recursion coefficients are

$$a_n = \frac{n(1+\beta) + 2\mu\beta}{1-\beta}, \qquad b_n = -\frac{\sqrt{\beta}}{1-\beta}\sqrt{(n+1)(n+2\mu)}. \tag{A5}$$

The discrete weight function reads as follows

$$\xi_k = (1-\beta)^{2\mu}(2\mu)_k \frac{\beta^k}{k!}. \tag{A6}$$

*The Krawtchouk polynomial:*

$$K_n^M(k;\gamma) = \sqrt{\frac{M!}{n!(M-n)!}} \left(\frac{\gamma}{1-\gamma}\right)^{n/2} \,_2F_1\left(\begin{matrix}-n,-k\\ -M\end{matrix}\bigg| \gamma^{-1}\right), \tag{A7}$$

where $1 > \gamma > 0$ and $k, n = 0, 1, 2, .., M$. The recursion coefficients are

$$a_n = M\gamma + n(1 - 2\gamma), \qquad b_n = -\sqrt{(n+1)(M-n)\gamma(1-\gamma)}. \tag{A8}$$

The discrete weight function reads as follows

$$\xi_k = (1-\gamma)^{N-k} \frac{\Gamma(M+1)\gamma^k}{\Gamma(M-k+1)k!}. \tag{A9}$$

## A.2 Discrete and continuous weight function mix

In this section of the Appendix, we give the essential properties of two orthogonal polynomials that have either a pure continuous measure (if $\mu > 0$) or a mix of continuous and discrete measures (if $\mu < 0$). The discrete spectrum points for both are $x_k^2 = -(k+\mu)^2$, which is finite and of size equal to $\lfloor -\mu \rfloor$ where $\lfloor z \rfloor$ is the largest integer less than $z$.

***The continuous dual Hahn polynomial:***

$$S_n^\mu(x^2; \alpha, \beta) = \sqrt{\frac{(\mu+\alpha)_n (\mu+\beta)_n}{n! (\alpha+\beta)_n}} \, {}_3F_2\left(\begin{matrix}-n, \mu+ix, \mu-ix \\ \mu+\alpha, \mu+\beta\end{matrix}\middle|1\right). \tag{A10}$$

It is a polynomial in $x^2$ of degree $n$ with $x_- = 0$ and $x_+ = +\infty$. If $\mu > 0$ then $\alpha$ and $\beta$ must be positive or a pair of complex conjugates with positive real parts. However, if $\mu < 0$ then $\alpha + \mu$ and $\beta + \mu$ should be positive or a pair of complex conjugates with positive real parts. The recursion coefficients are

$$a_n = (n+\mu+\alpha)(n+\mu+\beta) + n(n+\alpha+\beta-1) - \mu^2, \tag{A11a}$$

$$b_n = -\sqrt{(n+1)(n+\alpha+\beta)(n+\mu+\alpha)(n+\mu+\beta)}. \tag{A11b}$$

The continuous component of the weight function reads as follows

$$\sigma(x) = \frac{1}{2\pi} \frac{|\Gamma(\mu+ix)\Gamma(\alpha+ix)\Gamma(\beta+ix)/\Gamma(2ix)|^2}{\Gamma(\mu+\alpha)\Gamma(\mu+\beta)\Gamma(\alpha+\beta)}. \tag{A12}$$

The discrete component of the weight function is

$$\xi_k = 2\frac{\Gamma(\alpha-\mu)\Gamma(\beta-\mu)}{\Gamma(\alpha+\beta)\Gamma(1-2\mu)} \frac{(-\mu-k)}{(-1)^k k!} \frac{(\mu+\alpha)_k (\mu+\beta)_k (2\mu)_k}{(\mu-\alpha+1)_k (\mu-\beta+1)_k}, \tag{A13}$$

where $k = 0, 1, 2, ..., \lfloor -\mu \rfloor$.

***The Wilson polynomial:***

$$W_n^\mu(x^2; \nu, \alpha, \beta) =$$
$$= \sqrt{\left(\frac{2n+\mu+\nu+\alpha+\beta-1}{n+\mu+\nu+\alpha+\beta-1}\right) \frac{(\mu+\alpha)_n (\mu+\beta)_n (\mu+\nu)_n (\mu+\nu+\alpha+\beta)_n}{(\nu+\alpha)_n (\nu+\beta)_n (\alpha+\beta)_n n!}} \, {}_4F_3\left(\begin{matrix}-n, n+\mu+\nu+\alpha+\beta-1, \mu+ix, \mu-ix \\ \mu+\nu, \mu+\alpha, \mu+\beta\end{matrix}\middle|1\right) \tag{A14}$$

It is a polynomial in $x^2$ of degree $n$ with $x_- = 0$ and $x_+ = +\infty$. If $\mu > 0$ then $\{\nu, \alpha, \beta\}$ should be positive or complex conjugates with positive real parts. However, if $\mu < 0$ then $\{\nu, \alpha, \beta\} + \mu$ must be positive or complex conjugates with positive real parts. The recursion coefficients are

$$a_n = \frac{(n+\mu+\nu)(n+\mu+\alpha)(n+\mu+\beta)(n+\mu+\nu+\alpha+\beta-1)}{(2n+\mu+\nu+\alpha+\beta)(2n+\mu+\nu+\alpha+\beta-1)} + \frac{n(n+\nu+\alpha-1)(n+\nu+b-1)(n+\alpha+b-1)}{(2n+\mu+\nu+\alpha+\beta-1)(2n+\mu+\nu+\alpha+\beta-2)} - \mu^2, \quad \text{(A15a)}$$

$$b_n = -\frac{1}{2n+\mu+\nu+\alpha+\beta}\sqrt{\frac{(n+1)(n+\mu+\nu)(n+\alpha+\beta)(n+\mu+\alpha)(n+\mu+\beta)(n+\nu+\alpha)(n+\nu+\beta)(n+\mu+\nu+\alpha+\beta-1)}{(2n+\mu+\nu+\alpha+\beta-1)(2n+\mu+\nu+\alpha+\beta+1)}}. \quad \text{(A15b)}$$

The continuous component of the weight function reads as follows

$$\sigma(x) = \frac{1}{2\pi}\frac{\Gamma(\mu+\nu+\alpha+\beta)\left|\Gamma(\mu+\mathrm{i}x)\Gamma(\nu+\mathrm{i}x)\Gamma(\alpha+\mathrm{i}x)\Gamma(\beta+\mathrm{i}x)/\Gamma(2\mathrm{i}x)\right|^2}{\Gamma(\mu+\nu)\Gamma(\alpha+\beta)\Gamma(\mu+\alpha)\Gamma(\mu+\beta)\Gamma(\nu+\alpha)\Gamma(\nu+\beta)}. \quad \text{(A16)}$$

The discrete component of the weight function is

$$\xi_k = 2\frac{\Gamma(\mu+\nu+\alpha+\beta)\Gamma(\nu-\mu)\Gamma(\alpha-\mu)\Gamma(\beta-\mu)}{\Gamma(-2\mu+1)\Gamma(\alpha+\beta)\Gamma(\alpha+\nu)\Gamma(\beta+\nu)}\frac{(-\mu-k)(2\mu)_k(\mu+\nu)_k(\mu+\alpha)_k(\mu+\beta)_k}{(\mu-\nu+1)_k(\mu-\alpha+1)_k(\mu-\beta+1)_k k!}, \quad \text{(A17)}$$

where $k = 0, 1, 2, ..., \lfloor -\mu \rfloor$.

Therefore, if we approximate the function $f(x^2)$ by a Taylor series up to $(x^2)^l$ then $f(x^2)$ becomes a polynomial in $x^2$ of degree $l$, which we call $F_l(x^2)$. Thus, we can write

$$\int \rho(x) F_l(x^2) p_n(x^2) p_m(x^2) dx = \int_{x_-}^{x_+} \sigma(x) F_l(x^2) p_n(x^2) p_m(x^2) dx$$

$$+ \sum_{k=0}^{M} \xi_k F_l(x_k^2) p_n(x_k^2) p_m(x_k^2) = [F_l(J)]_{n,m}$$

Hence, the left side of Eq. (20) has the exact value as: $\lim_{l \to \infty} [F_l(J)]_{0,0} = [f(J)]_{0,0}$. To calculate the matrix $F_l(J)$, we proceed as follows. Let $\{\Lambda_{m,n}\}_{m=0}^{K-1}$ be the normalized eigenvector of the $K \times K$ finite submatrix of $J$ corresponding to the eigenvalue $\varepsilon_n$. Then, we can write $F_l(J) = \Lambda W \Lambda^T$ where $W$ is a diagonal matrix whose elements are: $W_{n,m} = \delta_{n,m} F_l(\varepsilon_n)$. For improved accuracy we take the matrix size $K$ as large as numerically possible.

## Table Caption

**Table 1:** Relative error in the evaluation of the infinite sum (14) for the function $f(x) = r^x / \Gamma(x+1)$ with $r = 3.0$ using Gauss quadrature associated with the Charlier and Meixner polynomials for several values of $N$. The polynomial parameter $\mu$ is fixed as $\mu = 2.0$ for both polynomials. However, we took several values for the Meixner polynomial parameter $\beta$. The relative error is calculated as $\left|\frac{e^r - \Sigma}{e^r + \Sigma}\right|$, where $\Sigma$ is the Gauss quadrature sum. For numerical computations, we used Mathcad® software version 14.0 with single precision.

**Table 2:** Relative error in the calculation of the finite sum (14) for the function $f(x) = (x+1) r^{x+1} / \Gamma(x+r+2)$ with $r = 3.0$ using Gauss quadrature associated with the Krawtchouk polynomial for $M = 100$ and for several values of $N$. We took several values for the polynomial parameter $\gamma$. The relative error is calculated as $\left|\frac{\text{Exc} - \Sigma}{\text{Exc} + \Sigma}\right|$, where Exc stands for the exact value obtained by formula (15) and $\Sigma$ is the Gauss quadrature sum.

**Table 3:** Relative error in the calculation of the integral-sum combination on the left side of Eq. (20) for the function $f(x) = x^m e^{-x/2}$ with $m = 3$. We used QG associated with the continuous dual Hahn polynomial $S_n^\mu(x^2; \alpha, \alpha)$ for $\mu = -3.5$ and for several values of the parameter $\alpha$. The relative error is calculated as $\left|\frac{\text{Exc} - \Sigma}{\text{Exc} + \Sigma}\right|$ where $\Sigma$ is the Gauss quadrature and Exc is $[f(J)]_{00}$ with a matrix size of 200 for $J$.

**Table 1**

| GQ | $N=2$ | $N=4$ | $N=7$ | $N=10$ | $N=15$ |
|---|---|---|---|---|---|
| Charlier | $5.694\times10^{-3}$ | $6.525\times10^{-6}$ | $4.165\times10^{-11}$ | $2.653\times10^{-16}$ | $8.844\times10^{-17}$ |
| Meixner ($\beta=0.2$) | $6.943\times10^{-3}$ | $1.231\times10^{-4}$ | $1.964\times10^{-7}$ | $1.522\times10^{-10}$ | $1.946\times10^{-15}$ |
| Meixner ($\beta=0.4$) | $3.900\times10^{-2}$ | $2.272\times10^{-3}$ | $3.192\times10^{-5}$ | $8.121\times10^{-7}$ | $1.1969\times10^{-9}$ |
| Meixner ($\beta=0.6$) | $9.541\times10^{-2}$ | $5.266\times10^{-3}$ | $1.131\times10^{-3}$ | $2.588\times10^{-5}$ | $8.008\times10^{-6}$ |

**Table 2**

| $\gamma$ | $N=10$ | $N=20$ | $N=30$ | $N=40$ | $N=50$ |
|---|---|---|---|---|---|
| 0.01 | $4.002\times10^{-11}$ | $7.725\times10^{-13}$ | $9.770\times10^{-15}$ | $2.220\times10^{-16}$ | $5.329\times10^{-15}$ |
| 0.10 | $3.600\times10^{-2}$ | $8.826\times10^{-6}$ | $2.469\times10^{-11}$ | $6.222\times10^{-12}$ | $5.390\times10^{-13}$ |
| 0.20 | $8.514\times10^{-1}$ | $4.065\times10^{-2}$ | $1.075\times10^{-4}$ | $9.438\times10^{-9}$ | $1.799\times10^{-14}$ |
| 0.30 | $9.999\times10^{-1}$ | $6.666\times10^{-1}$ | $4.314\times10^{-2}$ | $2.807\times10^{-4}$ | $8.968\times10^{-8}$ |

**Table 3**

| $\alpha+\mu$ | $N=10$ | $N=20$ | $N=30$ | $N=50$ | $N=100$ |
|---|---|---|---|---|---|
| 1.0 | $6.752\times10^{-5}$ | $4.338\times10^{-7}$ | $1.169\times10^{-8}$ | $6.258\times10^{-11}$ | $3.594\times10^{-12}$ |
| 2.0 | $2.012\times10^{-3}$ | $2.577\times10^{-5}$ | $9.999\times10^{-7}$ | $7.667\times10^{-9}$ | $2.048\times10^{-12}$ |
| 3.0 | $1.713\times10^{-2}$ | $4.119\times10^{-4}$ | $2.255\times10^{-5}$ | $2.584\times10^{-7}$ | $1.289\times10^{-10}$ |
| 4.0 | $7.529\times10^{-2}$ | $3.168\times10^{-3}$ | $2.403\times10^{-4}$ | $4.043\times10^{-6}$ | $3.385\times10^{-9}$ |
| 5.0 | $2.197\times10^{-1}$ | $1.494\times10^{-2}$ | $1.539\times10^{-3}$ | $3.743\times10^{-5}$ | $4.938\times10^{-8}$ |